\documentclass[a4paper,10pt]{amsart}
\usepackage{mathrsfs}
\usepackage[arrow,matrix]{xy}
\usepackage{amsmath,amssymb, amscd, amsthm,mathrsfs}
\usepackage[colorlinks=true,citecolor=blue,linkcolor=blue]{hyperref}
\theoremstyle{plain}
\newtheorem{thm}{Theorem}[section]
\newtheorem{cor}[thm]{Corollary}
\newtheorem{lem}[thm]{Lemma}
\newtheorem{prop}[thm]{Proposition}
\newtheorem{defn}[thm]{Definition}
\newtheorem{exa}[thm]{Example}

\newtheorem{rem}[thm]{Remark}
   \def\op{\oplus} \def\ot{\otimes}
\def\Hom{\operatorname {Hom}}\def\End{\operatorname {End}}\def\D{\mathcal{D}}
\def\RHom{\operatorname {RHom}}
\def\Ext{\operatorname {Ext}}

\def\k{\mathbf k}\def\M{\mathcal{M}}\def\Rat{\text{\rm Rat}}

\newcommand{\darrow}[2]{%
\genfrac{}{}{0pt}{}{\raisebox{-0.9ex}{$\xrightarrow{#1}$}}{\raisebox{0.9ex}{$\xleftarrow[#2]{}$}}}

\begin{document}
\title[Dualities of artinian coalgebras]{Dualities of
artinian coalgebras with applications to noetherian complete algebras}
\author{J.-W. He, B. Torrecillas, F. Van Oystaeyen and Y. Zhang}
\address{J.-W. He\newline \indent Department of Mathematics, Shaoxing College of Arts and Sciences, Shaoxing Zhejiang 312000,
China\newline \indent Department of Mathematics and Computer
Science, University of Antwerp, Middelheimlaan 1, B-2020 Antwerp,
Belgium} \email{jwhe@usx.edu.cn}
\address{B. Torrecillas\newline\indent Department of Algebra and
Analysis, University of Almeria, E-04071, Almeria, Spain}
\email{btorreci@ual.es}
\address{F. Van Oystaeyen\newline\indent Department of Mathematics and Computer
Science, University of Antwerp, Middelheimlaan 1, B-2020 Antwerp,
Belgium} \email{fred.vanoystaeyen@ua.ac.be}
\address{Y. Zhang\newline
\indent Department WNI, University of Hasselt, Universitaire Campus,
3590 Diepenbeek, Belgium} \email{yinhuo.zhang@uhasselt.be}

%\thanks{}
\date{}

\begin{abstract} A duality theorem of the bounded derived category
of quasi-finite comodules over an artinian coalgebra is established.
Let $A$ be a noetherian complete basic semiperfect algebra over an
algebraically closed field, and $C$ be its dual coalgebra. If $A$ is
Artin-Schelter regular, then the local cohomology of $A$ is
isomorphic to a shift of twisted bimodule ${}_1C_{\sigma^*}$ with
$\sigma$ a coalgebra automorphism. This yields that the balanced
dualinzing complex of $A$ is a shift of the twisted bimodule
${}_{\sigma^*}A_1$. If $\sigma$ is an inner automorphism, then $A$ is Calabi-Yau.
\end{abstract}

\keywords{Calabi-Yau (co)algebra, Duality, Derived category,
Artin-Schelter regular (co)algebra} \subjclass[2000]{16W30, 18E30,
18G10, 16E30}

\maketitle

\section*{Introduction}

The noncommutative dualizing complex, introduced by Yekutieli in
\cite{Yek1}, provides a powerful tool to study noncommutative algebras.
In order to determine a dualizing complex over a graded algebra, Yekutieli
introduced the concept of a {\it balanced} dualizing complex. However, the term `balanced' makes no sense for a general non-graded algebra. As an alternative, Van den Bergh introduced the concept of a {\it rigid}
dualizing complex (cf. \cite{VdB}) for a non-graded algebra. The
existence of a rigid (or balanced) dualizing complex of an algebra is related to the (twisted)
Calabi-Yau property of certain triangulated category (cf.
\cite{Gin}). Van den Bergh's results were generalized to noetherian
complete semilocal algebras (cf. \cite{CWZ,WZ}). Many good
properties of noetherian complete semilocal algebras were discovered through dualizing complexes.

However, any noetherian complete algebra $A$ with cofinite Jacobson
radical is the dual algebra of an artinian coalgebra $C$ (cf.
\cite{HR}). There exist certain duality properties between the
category of $A$-modules and the category of $C$-comodules. This
motivates us to study the balanced dualizing
complex, the Calabi-Yau property and the local cohomology of a
noetherian complete algebra through artinian coalgebras. To this aim, we first have to discuss some homological properties of
artinian coalgebras.

Let $C$ be an artinian coalgebra. In Section 1, we establish some dualities between triangulated
subcategories of derived categories of $C^*$-modules and those of
$C$-comodules. In particular, it turns out that the dual algebra
$C^*$ is Calabi-Yau if and only if $C$ is. Through the dualities
obtained in Section 1, we deduce the following duality theorem
of the bounded derived categories of left and of right $C$-comodules
when $C$ satisfies certain additional conditions in Section 2.

\noindent{\bf Theorem.} {\it Let $C$ be an artinian coalgebra. If
the following conditions are satisfied:
\begin{itemize}
  \item [(i)] the functors $\Gamma:{}_{C^*}\M\longrightarrow \M^C$
  and $\Gamma^\circ:\M_{C^*}\longrightarrow {}^C\M$ have finite
  cohomological dimensions;
  \item [(ii)] the coalgebra $C$ satisfies the left and the right
  $\chi$-condition,
\end{itemize}
then the functors $F=R\Gamma\circ(\ )^*$ and $G=R\Gamma^\circ\circ(\
)^*$ are dualities of triangulated categories:
$$D^b_{qf}({}^C\!\M)\darrow{\quad F\quad}{\quad G\quad}D^b_{qf}(M^C).$$}
In the theorem, ${}_{C^*}\M$ is the category of left $C^*$-modules, ${}^C\M$ is the category of left $C$-comodules, and $\mathcal{D}^b_{qf}({}^C\M)$ is the derived category of bounded complexes of left $C$-comodules with quasi-finite cohomology comodules. The $\chi$-condition in the theorem will be explained in the section
2, and $\Gamma=\Rat$ is the rational functor.

In Section 3 and Section 4, we focus on a class of artinian coalgebras
satisfying the $\chi$-condition, namely the Artin-Schelter regular
coalgebras. Let $A$ be a noetherian complete basic algebra over an
algebraically closed field, and let $C=A^\circ$ be its dual
coalgebra. If $A$ has cofinite Jacobson radical, then $A\cong C^*$.
If $A$ is Artin-Schelter regular, then $C$ is an Artin-Schelter regular
coalgebra. We have the following theorem (cf. Theorem
\ref{cthm3} and Corollary \ref{ccor3}).

\noindent{\bf Theorem.} {\it  Let $A$ be a noetherian complete basic
algebra with cofinite Jacobson radical over an algebraically closed
field, and $C$ be its dual coalgebra. Assume that $A$ is
Artin-Schelter regular of global dimension $n$. Then

{\rm(i)} there is a coalgebra automorphism $\sigma\in Aut(C)$ such
that $R\Gamma(A)\cong {}_1C_{\sigma^*}[-n]$ in $\D^b({}_{A}\M_{A})$;

{\rm(ii)} for any finitely generated left (or right) $A$-module,
$\dim\Ext^n_{A}(M, A)<\infty$; moreover, as vector spaces
$\Ext^n_{A}(M,A)\cong \Rat(M)^*$;

{\rm(iii)} for $i<n$, $\Ext_{A}^i(M,A)\cong\Ext^i_{A}(M/\Rat(M),A)$.

{\rm(iv)} if the automorphism $\sigma$ in (i) is inner, then $A$ is
Calabi-Yau.}

The items (ii) and (iii) can be viewed as generalizations of
\cite[Prop. 2.46(ii,iii)]{ATV} and \cite[Theorem 0.3(4)]{Zh}. The
item (i) says that $A$ has a balanced dualizing complex
${}_{\alpha}A_{1}[n]$ ($\alpha$ is an algebra automorphism), which
is similar to connected graded Artin-Schelter regular algebras (cf.
\cite{Yek1}).

Throughout $\k$ is an algebraically closed field of characteristic
zero. All the algebras and coalgebras involved are over $\k$;
unadorned $\ot$ means $\ot_\k$ and Hom means Hom$_\k$. Let $C$ be a
coalgebra, and $M$ and $N$ be left (or right) $C$-comodules. We use
$\Hom_C(M,N)$ to denote the set of left (or
right) $C$-comodule morphisms. To avoid
possible confusion, sometimes we use $\Hom_{C^{op}}(M,N)$ to
denote the set of right $C$-comodule morphisms for two $C$-bicomodules $M$ and $N$. We use similar Notations for modules over an algebra.

\section{Dualities between an artinian coalgebra and its dual algebra}
Let $C$ be an artinian coalgebra. Let $\Rat(\M)_{C^*}$ be the subcategory of $\M_{C^*}$ consisting of
rational modules. Then the abelian category $\Rat(\M)_{C^*}$ is equivalent
to the abelian category ${}^C\!\M$.

For a coalgebra $C$, ${}^CC$ is an injective object in ${}^C\!\M$,
or equivalently, $C_{C^*}$ is an injective object in
$\Rat(\M)_{C^*}$. In general, $C_{C^*}$ is not injective in
$\M_{C^*}$. However, we have the following property (see \cite[9.4]{BW} and \cite[Theorem 3.2]{CNV}).

\begin{prop} \label{prop1} The following are equivalent.
\begin{itemize}
  \item [(i)] $C_{C^*}$ is injective in $\M_{C^*}$;
  \item [(ii)] $C_{C^*}$ is an injective cogenerator of
  $\M_{C^*}$;
  \item [(iii)] ${}_{C^*}C$ is artinian;
  \item [(iv)] $C^*$ is left noetherian;
  \item [(v)] The injective hull of a rational left $C^*$-module is rational.
\end{itemize}
\end{prop}

From now on, unless stated otherwise, $C$ is both a left and a right artinian coalgebra.

Let $R$ and $S$ be two rings, ${}_RE_S$ an $R$-$S$-bimodule. Recall from \cite{AF} that an $R$-module
${}_RM$ (or an $S$-module $N_S$) is called an $E$-reflexive if the
natural morphism $$M\longrightarrow \Hom_S(\Hom_R(M,E),E)\qquad (\text{or},
N\longrightarrow \Hom_R(\Hom_S(N,E),E))$$ is an isomorphism. An
$R$-$S$-bimodule ${}_RE_S$ defines a {\it Morita duality} \cite[Sect. 24]{AF}
if
\begin{itemize}
\item [(i)] both ${}_RR$ and $S_S$ are $E$-reflexive;
\item [(ii)] every submodule and every quotient module of an
$E$-reflexive module is $E$-reflexive.
\end{itemize}

For a coalgebra $C$, consider the $C^*$-bimodule ${}_{C^*}C_{C^*}$.
Clearly we have $\End(C_{C^*})\cong C^*$ and
$\End({}_{C^*}C)\cong{C^*}^{op}$. Moreover, if $C$ is (both left and
right) artinian, then both ${}_{C^*}C$ and $C_{C^*}$ are injective
cogenerators by Prop. \ref{prop1}. So, the $C^*$-bimodule
${}_{C^*}C_{C^*}$ defines a Morita duality \cite[Theorem 24.1]{AF}.
Let $\mathcal{U}_{C^*}$ (resp. ${}_{C^*}\mathcal{U}$) be the
subcategory of $\M_{C^*}$ (resp. ${}_{C^*}\M$) consisting of
${}_{C^*}C_{C^*}$-reflexive modules. Then we have a duality:
\begin{equation}\label{ctag1}
\mathcal{U}_{C^*}\darrow{\Hom_{C^*}(-,C)}{\Hom_{C^*}(-,C)}{}_{C^*}\mathcal{U}.
\end{equation}
Let ${}_{C^*}\M_{fg}$ be the category of all finitely generated
left $C^*$-modules, and let ${}^C\M_{qf}$ be the category of all
 quasi-finite left $C$-comodule. From the duality (\ref{ctag1}), we
immediately obtain the following.

\begin{prop}\label{cprop1} Let $C$ be an artinian coalgebra. There is a duality
between abelian categories:
$${}_{C^*}\M_{fg}\darrow{\Hom_{C^*}(-,C)}{\quad(\ \ )^*\qquad}{}^C\M_{qf}.$$
\end{prop}

Since the dual algebra $C^*$ of an artinian coalgebra is noetherian
and semiperfect, the proposition above immediately follows the next
corollary (see also \cite[Prop. 3.6]{NNT} and \cite[Prop.
3.4]{CNV}).

\begin{cor} If $C$ is an artinian coalgebra, then gl.dim$C$=gl.dim$(C^*)$.
\end{cor}

Let $C$ be an arbitrary coalgebra, ${}_{C^*}M$ a $C^*$-module. There is
a right $C^*$-module morphism: $$\varphi_M:\Hom_{C^*}(M,C)\to M^*,\
f\mapsto \varepsilon\circ f.$$ If $M$ is a finitely generated
$C^*$-module, then the right $C^*$-module $\Hom_{C^*}(M,C)$ is in
fact a rational $C^*$-module. Moreover, if we view $\Hom_{C^*}(M,C)$
as a left $C$-comodule, then it is a quasi-finite comodule. Hence
the image of $\varphi_M$ is contained in $\Rat(M^*)$. Therefore, we
obtain a natural transformation:
\begin{equation}\label{ctag2}
    \varphi:\Hom_{C^*}(-,C)\longrightarrow \Rat\circ(\ )^*,
\end{equation}
of functors from ${}_{C^*}\M_{fg}$ to ${}^C\M_{qf}$.

\begin{lem}\label{clem1} If $C$ is an artinian coalgebra, then the natural transformation
$\varphi$ above is a natural isomorphism.
\end{lem}
\proof If ${}_{C^*}M$ is a finitely generated free module, then
$\varphi_M:\Hom_{C^*}(M,C)\longrightarrow \Rat\circ(M)^*$ is clearly
a isomorphism. For a general finitely generated module $M$, $M$ is
finitely presented since $C^*$ is noetherian:
$$\bigoplus_{\text{finite}} C^*\longrightarrow\bigoplus_{\text{finite}} C^*\longrightarrow M\longrightarrow0.$$
The statement follows from the following commutative diagram:
$$\xymatrix{
  0  \ar[r] & \Hom_{C^*}(M,C) \ar[d]^{\varphi_M}\ar[r]&\Hom_{C^*}(\displaystyle\bigoplus_{\text{finite}} C^*,C)\ar[d]^\varphi\ar[r]&
  \Hom_{C^*}(\displaystyle\bigoplus_{\text{finite}} C^*,C)\ar[d]^\varphi\\
  0 \ar[r] & \Rat\circ(M)^*\ar[r] &\Rat\circ\left(\displaystyle\bigoplus_\text{finite} C^*\right)^*\ar[r]&
  \Rat\circ\left(\displaystyle\bigoplus_\text{finite} C^*\right)^*. \qed}$$

Let $A$ be a noetherian algebra with Jacobson radical $J$ such that
$A_0=A/J$ is finite dimensional. Let ${}_AM$ be an $A$-module. An
element $m\in M$ is called a {\it torsion element} if $J^nm=0$ for
$n\gg0$. Let $\Gamma(M)=\{m\in M|m\ \text{is a torsion element}\}$.
Then $\Gamma(M)$ is a submodule of $M$. In fact, we have an additive
functor \cite{WZ}
$$\Gamma:{}_A\M\longrightarrow {}_A\M,$$ by sending an $A$-module to
its maximal torsion submodule. Clearly, $\Gamma$ is a left exact
functor. The functor $\Gamma$ has another representation
$\Gamma(M)=\underrightarrow{\lim}\Hom_A(A/J^n,M)$. We use $\Gamma^\circ$ to denote the torsion functor on the category of right $A$-modules.

Now let $C$ be a coalgebra, and let
$$\Rat:{}_{C^*}\M\longrightarrow {}_{C^*}\M$$ be the rational
functor. If $C$ is artinian, then by \cite[Prop. 3.1.1 and Remarks
3.1.2]{HR} every finite dimensional $C^*$-module is rational. Hence,
for a left $C^*$-module $M$, we have that $\Rat(M)$ is the sum of
all the finite dimensional submodules of $M$. On the other hand,
since the Jacobson radical $J=C_0^\bot$ and $C_0$ is finite
dimensional, $\Gamma(M)$ is also the sum of all the finite
dimensional submodules of $M$. Hence $\Gamma(M)\cong\Rat(M)$. So,
the functor $\Gamma$ is naturally isomorphic to the rational functor
$\Rat$. In what follow, we identify the right derived functor
$R\Gamma$ with $R\Rat$.

Let $C$ be an artinian coalgebra. Then $Soc(C)$ is finite
dimensional. This means that there are only finitely many non-isomorphic
simple right (or left) $C$-comodules. If ${}^CM$ is quasi-finite
then $Soc(M)$ is finite dimensional. Thus ${}^CM$ is finitely
cogenerated. This implies that ${}^C\M_{qf}$ is a thick subcategory
of ${}^C\M$. Hence $\D^+_{qf}({}^C\M)$, the derived category of
bounded below complexes of left $C$-comodule with quasi-finite
cohomology comodules, is a full triangulated subcategory of
$\D^+({}^C\M)$. Also, since $C^*$ is noetherian,
$\D^-_{fg}({}_{C^*}\M)$, the derived category of bounded above
complexes of left $C^*$-modules with finitely generated cohomology
modules, is a full triangulated subcategory of $\D^-({}_{C^*}\M)$.
The duality in Prop. \ref{cprop1} induces a duality of derived
categories.

\begin{prop}\label{cprop2} Let $C$ be an artinian coalgebra. We have dualities of
triangulated categories:
$$\D^-_{fg}({}_{C^*}\M)\darrow{\ R\Gamma^\circ\circ(\ \ )^*\ }{\quad(\ \
)^*\quad}\D^+_{qf}({}^C\M),\quad \D^b_{fg}({}_{C^*}\M)\darrow{\
R\Gamma^\circ\circ(\ \ )^*\ }{\quad(\ \ )^*\quad}\D^b_{qf}({}^C\M).$$
\end{prop}
\proof Since $C$ is artinian, $\D^+_{qf}({}^C\M)$ is equivalent to
$\D^+({}^C\M_{qf})$, the derived category of complexes of
quasi-finite comodules. In the dual case, $C^*$ is noetherian, and
$\D^-_{fg}({}_{C^*}\M)$ is equivalent to $\D^-({}_{C^*}\M_{fg})$.
By Prop. \ref{cprop1}, we have the following duality
$$\D^-({}_{C^*}\M_{fg})\darrow{\Hom_{C^*}(-,C)}{\quad(\ \
)^*\qquad}\D^+({}^C\M_{qf}).$$ Now using Lemma \ref{clem1} one
may check without difficulty that the composition
$$\xymatrix@C=0.8cm{
  \D^-_{fg}({}_{C^*}\M) \ar[r]^\cong & \D^-({}_{C^*}\M_{fg})\ar[rr]^{\Hom_{C^*}(-,C)} &&\D^+({}^C\M_{qf})  \ar[r]^\cong &\D^+_{qf}({}^C\M)
  }$$ is naturally isomorphic to the functor $R\Gamma^\circ\circ(\ )^*$.
Moreover, one sees that $R\Gamma^\circ\circ(\ )^*$ sends bounded complexes to
bounded complexes. \qed

\begin{cor}\label{ccor1} Let $C$ be an artinian coalgebra. Then we have a
duality of triangulated categories: $$\D^b_{fd}({}_{C^*}\M)\darrow{\
R\Gamma\circ(\ \ )^*\ }{\quad(\ \ )^*\quad}\D^b_{fd}({}^C\M).$$
\end{cor}
\proof It suffices to show that
$R\Gamma\circ(M)^*\in\D^b_{fd}({}^C\M)$ for all finite dimensional
left $C^*$-module $M$. Since $C^*$ is noetherian and complete with
respect to the radical filtration, the Jacobson radical $J$ of $C^*$
satisfies Artin-Rees condition. Hence the injective envelop of a
$J$-torsion module is still $J$-torsion (cf. \cite[Theorem
3.2]{CNV}). Now $M^*$ is a $J$-torsion module. We have an injective
resolution of $M^*$ with each component being $J$-torsion. Hence
$R\Gamma(M^*)$ is quasi-isomorphic to $M^*$, that is,
$R\Gamma(M^*)\in\D^b_{fd}({}^C\M)$.\qed

Recall that a $\k$-linear category $\mathcal{C}$ is said to be {\it Hom-finite}, if for any $X,Y\in \mathcal{C}$, $\Hom_{\mathcal{T}}(X,Y)$ is a finite dimensional $\k$-vector space; a Hom-finite $\k$-linear triangulated category $\mathcal{T}$ is called a {\it Calabi-Yau category of dimension $n$} if, for any objects $X,Y\in\mathcal{T}$, there is a natural isomorphism $\Hom_{\mathcal{T}}(X,Y)\cong\Hom_{\mathcal{T}}(Y,X[n])^*$; an algebra $A$ is called a (left) {\it Calabi-Yau of dimension $n$} (simply, CY-$n$) if

(i) $\D^b_{fd}({}_A\M)$ is Hom-finite;

(ii) $\D^b_{fd}({}_A\M)$ is a Calabi-Yau category of dimension $n$.

Note that the CY property of an algebra is always left-right symmetric. Thus we simply say that an algebra is CY-$n$.

In the dual case, we say that a coalgebra is (left) {\it CY-$n$} if

(i) $\D^b_{fd}({}^C\M)$ is Hom-finite;

(ii) $\D^b_{fd}({}^C\M)$ is a Calabi-Yau category of dimension $n$.

In general, we don't know  whether the CY property of a coalgebra is left-right symmetric. But an artinian coalgebra is left CY if and only if it is right CY. In fact, from the definitions we have the following.

\begin{cor} \label{ccor2} Let $C$ be an artinian coalgebra. Then $C$ is (left) CY-$n$ if and only if $C^*$ is (left) CY-$n$.
\end{cor}

\begin{exa} \label{exa1} {\rm Let $C$ be the path coalgebra \cite{Chin} of the
quiver $Q$ with one vertex and one arrow. Then the dual algebra
$C^*$ is the formal power series algebra $\k[[x]]$. It is well
known that $\k[[x]]$ is a CY-1 algebra. Hence $C$ is CY-1
coalgebra.}
\end{exa}

\section{Dualities of comodules over an artinian coalgebra}

In this section, we establish a duality of the derived
categories of left $C$-comodules and of right $C$-comodules by using
the results obtained in Section 1.

As ${}^C\!\M_{qf}$ is a thick subcategory of ${}^C\!\M$,
$D^b_{qf}({}^C\!\M)$ is a full triangulated subcategory of
$D^b({}^C\!\M)$. Consider the following functors:
$$\xymatrix{
  D^b({}^C\M) \ar[r]^{(\ )^*} & D^b({}_{C^*}\M) \ar[r]^{R\Gamma} &
  D^+(\M^C).  }$$
We want to know when the composite $R\Gamma\circ (\ )^*$ has its
image in $D^b(\M^C)$. If this happens, does the
restriction of $R\Gamma\circ (\ )^*$ to the subcategory
$D^b_{qf}({}^C\M)$ result a functor $D^b_{qf}({}^C\M)\longrightarrow
D^b_{qf}(\M^C)$?

We need to translate some concepts relative to noncommutative
algebras to coalgebras.

\begin{defn} {\rm Let $C$ be an artinian coalgebra. We
say that a quasi-finite left $C$-comodule $M$ satisfies the {\it
$\chi$-condition} if, for every simple left $C$-comodule $S$,
$\Ext_C^i(M,S)$ is finite dimensional for all $i\ge0$. We say that a
coalgebra $C$ satisfies the {\it left $\chi$-condition} if every
quasi-finite left $C$-comodule satisfies the $\chi$-condition.
Similarly, we can define the {\it right $\chi$-condition}.}
\end{defn}

The $\chi$-condition on a coalgebra is dual to a similar condition
on a noetherian algebra, which is originally introduced in \cite{AZ}
for noetherian graded algebras and was extended to nongraded algebras in
\cite{WZ,CWZ}.

\begin{defn} {\rm Let $A$ be a noetherian algebra with Jacobson radical
$J$ such that $A_0=A/J$ is finite dimensional. A finitely generated
left $A$-module $M$ is said to satisfy {\it the $\chi$-condition} if
$\Ext_A^i(A_0,M)$ is finite dimensional for all $i\ge0$. $A$ is said
to satisfy the {\it left $\chi$-condition} if every finitely
generated left $A$-module satisfies the $\chi$-condition.}
\end{defn}

Let ${}_{C^*}N$ be a module. There is a natural $C^*$-module
morphism $$\eta_N:N\to \Gamma(N^*)^*,\ \eta(n)(f)=f(n)$$ for $n\in
N$ and $f\in\Gamma(N^*)$. Observe that if $N$ is a finitely generated free $C^*$-module, then $\eta_N$ is an isomorphism.

\begin{lem}\label{clem2} If ${}_{C^*}N$ is a finitely generated module, then
$R\Gamma(N^*)$ is quasi-isomorphic to $\Gamma(N^*)$.
\end{lem}
\proof Choose a projective resolution of $N$ as follows:
$$0\longleftarrow N\longleftarrow P^0\longleftarrow P^{-1}\longleftarrow\cdots\longleftarrow
P^{-i}\longleftarrow\cdots,$$ where each $P^{-i}$ is a finitely
generated free $C^*$-module. Consider the following commutative
diagram:
$$\xymatrix{
0  & N\ar[d]^{\eta_N}\ar[l]&P^0\ar[l]\ar[d]^{\eta_{P^0}}&\cdots\ar[d]\ar[l]&P^i\ar[d]^{\eta_{P^i}}\ar[l]&\cdots\ar[l]\\
  0 & \Gamma(N^*)^*\ar[l]& \Gamma({P^0}^*)^*\ar[l]&\cdots\ar[l]&\Gamma({P^{-i}}^*)^*\ar[l]&\cdots\ar[l] }$$
Since $P^{-i}$ is a finitely generated free $C^*$-module,
$\eta_{P^{-i}}$ is an isomorphism for all $i\ge0$. Hence $\eta_N$ is
an isomorphism since $\Gamma$ is left exact. Now exactness of the
top row implies that the bottom row is exact too. Then the sequence
$$\xymatrix{
  0\ar[r] & \Gamma(N^*)\ar[r]& \Gamma({P^0}^*)\ar[r]&\cdots\ar[r]&\Gamma({P^{-i}}^*)\ar[r]&\cdots
  }$$ is also exact. Note that the sequence ${P^\bullet}^*$ is an
  injective resolution of $N^*$. Hence $R\Gamma(N^*)$ is
  quasi-isomorphic to $\Gamma(N^*)$. \qed

\begin{prop} $C$ satisfies the left $\chi$-condition if and only if
$C^*$ satisfies the left $\chi$-condition.
\end{prop}
\proof Note that the left $\chi$-condition on $C$ is equivalent to
the condition that for any quasi-finite left $C$-comodule $M$ and
any finite dimensional comodule $K$, $\Ext^i_C(M,K)$ is finite
dimensional for all $i\ge0$. Let $N$ be any finitely generated
$C^*$-module and $L$ be a finite dimensional $C^*$-module. By Prop.
\ref{cprop2}, $\Ext^i_{C^*}(L,N)=\Hom_{\D^b_{fg}({}_{C^*}\M)}(L,N[i])\cong\Hom_{\D^b_{qf}({}^C\M)}(R\Gamma(N^*),R\Gamma(L^*)[i])$. Since $R\Gamma(N^*)$ and $R\Gamma(L^*)$ are quasi-isomorphic with
$N^*$ and $L^*$ respectively (see the proof of Corollary \ref{ccor1}), we have isomorphism:
$R\Gamma(L^*)$ is quasi-isomorphic to $L^*$,
$\Ext^i_{C^*}(L,N)\cong\Hom_{\D^b_{qf}({}^C\M)}(\Gamma(N^*),L^*[i])$.
Now if $C$ satisfies the left $\chi$-condition,
$\Hom_{\D^b_{qf}({}^C\M)}(\Gamma(N^*),L^*[i])$ is finite dimensional
for all $i\ge0$. It follows that $\Ext^i_{C^*}(L,N)$ is finite dimensional for
all $i\ge0$. Similarly we see the converse is also true.\qed

Now we may establish a Morita-type duality of the derived categories of
comodules. A Morita-type equivalence of the derived categories of
comodules was shown in \cite{Fa}.

\begin{thm}\label{cthm2} Let $C$ be an artinian coalgebra. If the following conditions are satisfied:
\begin{itemize}
  \item [(i)] the functors $\Gamma:{}_{C^*}\M\longrightarrow \M^C$
  and $\Gamma^\circ:\M_{C^*}\longrightarrow {}^C\M$ have finite
  cohomological dimensions;
  \item [(ii)] the coalgebra $C$ satisfies the left and the right
  $\chi$-conditions,
\end{itemize}
then the functors $F=R\Gamma\circ(\ )^*$ and $G=R\Gamma^\circ\circ(\
)^*$ are dualities of triangulated categories:
$$D^b_{qf}({}^C\!\M)\darrow{\quad F\quad}{\quad G\quad}D^b_{qf}(M^C).$$
\end{thm}
\proof First of all, we have to show that the functors $F$ and $G$
are well-defined. For an object $X\in D^b_{qf}({}^C\M)$, $X^*$ lies in
$D^b_{fg}({}_{C^*}\M)$. By (i), $\Gamma$ has finite cohomological
dimension. Then $R\Gamma(X^*)$ is in $D^b(\M^C)$. We show
that $R\Gamma(X^*)$, in fact, belongs to $D_{qf}^b(\M^C)$. Since $X\in
D^b_{qf}({}^C\!\M)$, there is no harm to assume that $X$ is a
finitely cogenerated $C$-comodule. Now $X^*$ is a finitely generated
$C^*$-module. Let $I^0\longrightarrow I^1\longrightarrow
I^2\longrightarrow\cdots\longrightarrow I^i\longrightarrow\cdots$ be
a minimal injective resolution of $X^*$. By (ii), $C$ satisfies the left
$\chi$-condition. Then $\Hom_{C^*}(S,I^i)$ is finite dimensional for
all simple left $C^*$-module $S$ and all $i\ge0$. Hence $Soc(I^i)$
is finite dimensional for all $i\ge0$. Note that
$Soc(I^i)\subseteq\Rat(I^i)=\Gamma(I^i)$ and $Soc(I^i)$ is essential
in $\Gamma(I^i)$. We obtain that $\Gamma(I^i)$ is finitely
cogenerated as a right $C$-module for all $i\ge0$. Then
$R\Gamma(X^*)$, which is quasi-isomorphic to $\Gamma(I^\bullet)$,
must have finitely cogenerated $C$-comodules as its cohomologies.
Hence $F(X)=R\Gamma(X^*)\in D^b_{qf}(\M^C)$. Similarly, we can show
that $G(Y)$ lies in $D^b_{qf}({}^C\M)$ for any $Y\in D^b_{qf}(\M^C)$.

By assumption, $C^*$ is noetherian, complete and semiperfect.
Thus $C^*/J$ is semisimple and finite dimensional. Applying
Theorem 4.1 of \cite{WZ} to our case, we obtain a
duality of the following triangulated categories:
$$\D^b_{fg}(\M_{C^*})\darrow{\Hom_{C^*}(R\Gamma^\circ(-),C)}{\Hom_{C^*}(R\Gamma(-),C)\ }\D^b_{fg}({}_{C^*}\M).$$
On the other hand, following Prop. \ref{cprop2}, we have dualities
$$\D^b_{fg}({}_{C^*}\M)\darrow{\ R\Gamma\circ(\ \ )^*\ }{\quad(\ \
)^*\quad}\D^b_{qf}({}^C\M),\ \ \D^b_{fg}(\M_{C^*})\darrow{\
R\Gamma^\circ\circ(\ \ )^*\ }{\ \quad(\ \ )^*\quad\
}\D^b_{qf}(\M^C).$$ We want to show that the composition
$$\xymatrix@C=1.2cm{
 \Psi: \D^b_{qf}({}^C\M)\ar[r]^{(\ \ )^*} & \D^b_{fg}({}_{C^*}\M) \ar[rr]^{\Hom_{C^*}(R\Gamma(-),C)}
  &&\D^b_{fg}(\M_{C^*})\ar[r]^{R\Gamma\circ(\ \ )^*}
  &\D^b_{qf}(\M^C)}$$
is naturally isomorphic to the functor $F=R\Gamma\circ(\ )^*$. For
any $X\in \D^b_{qf}({}^C\M)$, from the last paragraph we see
$R\Gamma(X^*)\in \D^b_{qf}(\M^C)$. We have natural isomorphisms
$$\Hom_{C^*}(R\Gamma(X^*),C)^*=\Hom_{C^*}(\Hom_{C^*}(R\Gamma(X^*),C),C)\cong
R\Gamma(X^*)$$ since ${}_{C^*}C_{C^*}$ defines a Morita duality.
Hence we have natural isomorphisms in $\D^b_{qf}(\M^C)$
$$R\Gamma(\Hom_{C^*}(R\Gamma(X^*),C)^*)\cong R\Gamma(R\Gamma(X^*))\cong
R\Gamma(X^*),$$ where the last isomorphism holds because
$R\Gamma(X^*)$ is a complex of $J$-torsion modules as an object in
$\D^b_{fg}({}_{C^*}\M)$. So, $\Psi$ and $F$ are naturally isomorphic.
Similarly, we see $G$ is naturally isomorphic to the following composite functor:
$$\xymatrix@C=1.2cm{
\Psi: \D^b_{qf}(\M^C)\ar[r]^{(\ \ )^*} & \D^b_{fg}(\M_{C^*})
\ar[rr]^{\Hom_{C^*}(R\Gamma^\circ(-),C)}
  &&\D^b_{fg}({}_{C^*}\M)\ar[r]^{R\Gamma^\circ\circ(\ \ )^*}
  &\D^b_{qf}({}^C\M)}.$$ Then we get the desired results. \qed

\begin{rem} When $C$ is not artinian, there maybe exist certain
dualities between the derived categories of comodules. For example,
let $C$ be a (both left and right) semiperfect coalgebra, that is,
the categories ${}^C\M$ and $\M^C$ have enough projective objects
(cf. \cite{Lin}). Then the rational functor $\Rat$ is exact
\cite{GN1,GN2,NT}, and $\Rat\circ(\ )^*$ gives a Colby-Fuller
duality \cite[Theorem 3.5]{GN1}:
$$\M^C\darrow{\quad\Rat\circ(\ )^*\quad}{\quad\Rat\circ(\
)^*\quad}{}^C\!\!\M.$$ Of course, this duality induces a duality
between certain triangulated subcategories of the derived categories
of $\M^C$ and ${}^C\M$ respectively.
\end{rem}

There are two questions arising from the above theorem: (i) when does an artinian coalgebra satisfy the $\chi$-conditon? (ii) when is $F({}^CC)$ isomorphic to $C^C$ in $\D^b_{qf}(\M^C)$? We deal with these questions in the rest of the paper.

\section{Coalgebras satisfy the $\chi$-condition}

In this section, we introduce a class of coalgebras dual to
Artin-Schelter algebra. These coalgebras satisfy the $\chi$-condition
of Theorem \ref{cthm2}.

The classical concept of Artin-Shelter (AS, for short) regular
algebra \cite{AZ} is defined over graded algebras. We may extend
this concept from graded algebras to semiperfect algebras. We say
that a noetherian semiperfect algebra $A$ is {\it left AS-regular}
(cf. \cite{CWZ}) if $A$ has finite global dimension $d$ and for
every left simple $A$-module $S$, one has
$$\Ext^i_A(S,A)=\left\{
                                              \begin{array}{ll}
                                                0, & \hbox{if $i\neq d$;} \\
                                                T, & \hbox{if $i=d$,}
                                              \end{array}
                                            \right.
$$ where $T$ is a right simple $A$-module. Similarly, one can define a {\it right AS-regular} algebra.

In the dual case, we have the following definition.

\begin{defn} {\rm Let $C$ be an artinian coalgebra with global dimension
$d<\infty$. We say that $C$ is {\it right AS regular} if, for every
simple left $C^*$-module $S$, we have
$$\Ext^i_{C^*}(C,S)=\left\{
                              \begin{array}{ll}
                                0, & \hbox{if $i\neq d$;} \\
                                T, & \hbox{if $i=d$,}
                              \end{array}
                            \right.
$$ where $T$ is a simple left $C^*$-module.

Similarly, we may define a {\it left AS-regular} coalgebra.}
\end{defn}

\begin{rem} In the definition above, we may use alternatively the extensions in the category of comodules to define an AS-regular coalgebra. In fact, we have $\Ext^i_C(C,S)=\Ext^i_{C^*}(C,S)$ since the subcategory
${}_{C^*}\Rat(\M)$ is a thick subcategory. Also $T=\Ext^d_C(C,S)$ in
the definition above is a right $C$-comodule.
\end{rem}

Following Prop. \ref{cprop1}, one can check without difficulty that the following holds.

\begin{prop} A coalgebra $C$ is left AS-regular if and
only if $C^*$ is left AS-regular.
\end{prop}

A graded AS-regular algebra must satisfy the $\chi$-condition (cf.
\cite{AZ}, Sec.8). Similarly, a nongraded AS-regular algebra also
satisfies the $\chi$-condition (cf. \cite{CWZ}). In the dual case, we have
the following.

\begin{prop}\label{prop5} If $C$ is a left AS-regular coalgebra, then $C$
satisfies the left $\chi$-condition.
\end{prop}
\proof Note that the AS-regularity of $C$ and that of $C^*$ are
equivalent. \qed

We will see that the concept of an AS-regular algebra (coalgebra) is
left-right symmetric. Let $A$ be a left AS-regular algebra of global
dimension $d$. If ${}_{A}S$ is a simple module, then we write
$S^\natural$ for the right simple module $\Ext^d_{A}(S,A)$.

\begin{lem}\label{lem2} Let $A$ be a left AS-regular algebra of global dimension $d$. Then for two left simple $A$-modules $S$ and $T$, $S^\natural\cong T^\natural$ if and only if $S\cong T$.
\end{lem}
\proof  Let
$$0\longrightarrow P^{-d}\longrightarrow\cdots\longrightarrow
P^{-1}\longrightarrow P^0\longrightarrow S\longrightarrow0,$$ and
$$0\longrightarrow Q^{-d}\longrightarrow\cdots\longrightarrow
Q^{-1}\longrightarrow Q^0\longrightarrow T\longrightarrow0$$ be
minimal projective resolutions of $S$ and $T$ respectively. Applying
the functor $\Hom_{A}(-,A)$ to both resolutions, we obtain exact
sequences:{\small
\begin{equation}\label{tag1}
 0\longrightarrow   \Hom_{A}(P^{0},A)\longrightarrow\Hom_{A}(P^{-1},A)\longrightarrow\cdots\longrightarrow
\Hom_{A}(P^{-d},A)\longrightarrow S^\natural\longrightarrow0,
\end{equation}} and {\small
\begin{equation}\label{tag2}
0\longrightarrow
\Hom_{A}(Q^{0},A)\longrightarrow\Hom_{A}(Q^{-1},A)\longrightarrow\cdots\longrightarrow
\Hom_{A}(Q^{-d},A)\longrightarrow T^\natural\longrightarrow0.
\end{equation}} Since $P^i$ and $Q^i$ are finitely generated projective modules
for all $i\ge0$, the above sequences are minimal projective
resolutions of $S^\natural$ and $T^\natural$ respectively. Suppose
$S^\natural\cong T^\natural$. Then the sequences (\ref{tag1}) and
(\ref{tag2}) are isomorphic. Notice that we have isomorphisms of
complexes
$$\Hom_{A}(\Hom_{A}(P^\bullet,A),A)\cong P^\bullet$$ and
$$\Hom_{A}(\Hom_{A}(Q^\bullet,A),A)\cong Q^\bullet.$$
Thus we
obtain $S\cong T$. \qed

\begin{prop}\label{prop7} A noetherian semiperfect algebra is
left AS-regular if and only if it is right AS-regular.

As a consequence, an artinian coalgebra is left AS-regular if and
only if it is right AS-regular.
\end{prop}
\proof Suppose that $A$ is left AS-regular. Let ${}_AS$ be a simple
module. By the proof of Lemma \ref{lem2}, the right simple module
$S^\natural$ has a minimal projective resolution (\ref{tag1}), and
$\Ext^i_{A}(S^\natural,A)=0$ for $i\neq d$ and
$\Ext^d_{A}(S^\natural,A)=S$. Since there are only finitely many
nonisomorphic simple left (and right) $A$-modules, by Lemma
\ref{lem2}, for every right simple $A$-module $K$ there is a left
simple $A$-module $S$ such that $K=S^\natural$. Hence $A$ is right
AS-regular. \qed

In view of the proposition above, we may omit the prefix ``left'' and
``right'' and just say an AS-regular (co)algebra.

Let $C$ be an artinian coalgebra satisfying left and right $\chi$-conditions. The
local cohomology of $C^*$ (relative to the Jacobson radical)
provides a duality between certain triangulated categories (see Theorem
\ref{cthm2}). Furthermore, it gives a `balanced' dualizing complex
of $C^*$ (cf. \cite{CWZ}, for the terminology). If $C$ is an
AS-regular coalgebra, we can compute the local cohomology of $C^*$.
Recall that a coalgebra is {\it basic} if the dual of any simple
subcoalgebra is a division algebra (cf. \cite{CM}). If $C$ is a basic
artinian coalgebra, then $C^*$ is a noetherian basic algebra.

\begin{thm}\label{cthm3} If $C$ is a basic AS-regular coalgebra, then
there is a coalgebra automorphism $\sigma\in Aut(C)$ and a
nonnegative integer $n$ such that $R\Gamma(C^*)\cong
{}_1C_{\sigma^*}[-n]$ in $\D^b({}_{C^*}\M_{C^*})$.
\end{thm}
\proof When $C^*$ is a local algebra, the result can be deduced from
\cite[Cor. 3.9]{CWZ}. But it seems that we could not extend the
proof in \cite{CWZ} directly to the general case.

Since we work over an algebraically closed field, any simple
comodule over the basic coalgebra $C$ is one-dimensional. Assume $C$
is of global dimension $n$. Let
\begin{equation}\label{ctag3}
    0\longrightarrow {}_{C^*}C^*\longrightarrow I^0\longrightarrow
I^1\longrightarrow\cdots\longrightarrow I^n\longrightarrow0
\end{equation}
be the minimal injective resolution of the left module
${}_{C^*}C^*$. Let ${}_{C^*}S$ be any simple $C^*$-module. We have
$\Ext^i_{C^*}(S,C^*)\cong\Hom_{C^*}(S,I^i)$ for all $i\ge0$. By the
AS-regularity of $C$, we obtain $\Hom_{C^*}(S,I^i)=0$ for
$i\neq n$ and $\dim\Hom_{C^*}(S,I^n)=1$. Then we have $Soc(I^i)=0$
for all $i<n$ and $Soc(I^n)={}_{C^*}C_0$, where $C_0$ is the
coradical of $C$. Hence $I^i$ is $J$-torsion free for all $i<d$, and
$I^n={}_{C^*}C\op \bar{I}^n$ for some $J$-torsion free module
$\bar{I}^n$. Therefore we get $R\Gamma(C^*)\cong {}_{C^*}C[-n]$ in
$\D^b({}_{C^*}\M)$. Note that $R\Gamma(C^*)$ is an object in
$\D^b({}_{C^*}\M_{C^*})$. Let $T=H^{-n}(R\Gamma(C^*)^*)$. Then $T$
is a $C^*$-bimodule, and $T$ is free as a right $C^*$-module. Since
$C$ (and dually $C^*$) satisfies left and right $\chi$-conditions, by \cite[Theorem
4.1]{WZ}, $T$ is a dualizing complex of $C^*$. Since $C^*$ is of
finite global dimension, $C^*$ itself is a dualizing complex. By
\cite[Theorem 4.5]{Yek}, $T$ is a tilting complex. Now applying the right
version of \cite[Prop. 2.3]{RZ}, we obtain an algebra automorphism
$\alpha\in Aut(C^*)$ such that $T\cong {}_\alpha C^*_1$. Since $C$
is reflexive, there is a unique coalgebra automorphism $\sigma\in
Aut(C)$ such that $\sigma^*=\alpha$. Let $U=H^n(R\Gamma(C^*))$. It is clear that $T=U^*$. Since $U\cong{}_{C^*}C$ as a left $C^*$-module,
$U^*\cong\Hom_{C^*}(U,C)$ as $C^*$-bimodules. Since
${}_{C^*}C_{C^*}$ defines a Morita duality, the canonical morphism
$U\longrightarrow\Hom_{{C^*}^{op}}(\Hom_{C^*}(U,C),C)$ is an
isomorphism of $C^*$-bimodules. Thus we have $U\cong
\Hom_{{C^*}^{op}}(T,C)\cong {}_1C_{\alpha}={}_1C_{\sigma^*}$ as
$C^*$-bimodules. Therefore $R\Gamma(C^*)\cong {}_1C_{\sigma^*}[-n]$.
\qed

\begin{rem}\label{crem1} {\rm(i)} We call the automorphism $\sigma\in Aut(C)$ {\rm the Nakayama
automorphism} of $C$, and call its dual automorphism $\sigma^*$ {\rm
the Nakayama automorphism} of $C^*$. Note that Nakayama automorphism
is unique up to inner automorphisms.

{\rm(ii)} The theorem implies that ${}_{\sigma^*}C^*_{\ 1}[n]$ is
the balanced dualizing complex (cf. \cite{CWZ}) of $C^*$.

{\rm(iii)} Let $A$ be a noetherian complete (with respect to the
Jacobson radical) basic algebra, and let $C(:=A^\circ)$ be its dual
coalgebra. Then $C$ is an artinian basic coalgebra and $A=C^*$ (cf.
\cite[Prop. 4.3.1]{HR}). So, the preceding theorem applies to all
the noetherian complete AS-regular basic algebra.
\end{rem}

From Theorem \ref{cthm3}, one can deduce the following finiteness properties of extension groups of finitely generated $C^*$-modules, which can be viewed as a generalization of \cite[Prop. 2.46(ii,iii)]{ATV} and \cite[Theorem 0.3(4)]{Zh}.

\begin{cor} Let $C$ be a basic AS-regular coalgebra of global
dimension $n$. If $M$ is a finitely generated left (or right)
$C^*$-module, then:

{\rm(i)} $\dim\Ext^n_{C^*}(M, C^*)<\infty$; moreover, as vector
spaces $\Ext^n_{C^*}(M,C^*)\cong \Rat(M)^*$;

{\rm(ii)} for $i<n$,
$\Ext_{C^*}^i(M,C^*)\cong\Ext^i_{C^*}(M/\Rat(M),C^*)$.
\end{cor}
\proof (i) By Theorem \ref{cthm3}, $R\Gamma(C^*)\cong
{}_{\sigma^*}C^*_1[n]$ in $\D^b({}_{C^*}\M_{C^*})$. Applying the
local duality theorem \cite[Prop 3.4]{CWZ}, we have
$$R\Gamma(M)^*\cong\RHom_{C^*}(M,{}_{\sigma^*}C^*_{\ 1}[n]).$$ Taking the 0-th cohomology on both sides of the complexes above, we obtain
$$\Rat(M)^*\cong\Hom_{\D^b({}_{C^*}\M)}(M,{}_{\sigma^*}C^*_{\ 1}[n])\cong\Ext^n_{C^*}(M,{}_{\sigma^*}C^*_{\
1}).$$ Since $M$ is finitely generated, $\Rat(M)$ is finite
dimensional.

(ii) From the proof of the preceding theorem, in the minimal
injective resolution (\ref{ctag3}) of $C^*$, $I^i$ is torsion free
for all $i<n$. Then we have
$\Hom_{C^*}(M,I^i)\cong\Hom_{C^*}(M/\Rat(M),I^i)$. Hence (ii)
follows. \qed

We end this section with an example of basic artinian AS-regular
coalgebra.

\begin{exa} {\rm Let $Q$ be the following quiver: $$
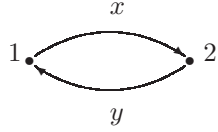
\begin{picture}(200,40)(0,0)
\put(100,10){\circle*{3}}\put(92,10){\makebox{$1$}}\put(160,10){\circle*{3}}\put(165,10){\makebox{$2$}}
\qbezier(102,12)(130,30)(158,12)\qbezier(158,8)(130,-10)(102,8)
\put(102,8){$\vector(-4,3){0}$}\put(158,12){$\vector(4,-3){0}$}
\put(130,28){\makebox{$x$}}\put(130,-12){\makebox{$y$}}
\end{picture}$$

Let $C$ be the path coalgebra $CQ$. Let $S_1$ and $S_2$ be the
simple left $C$-comodules corresponding to the vertices, and
$e_1,e_2\in C^*$ be the idempotents corresponding to the vertices.
Now the injective envelop of $S_1$ is $e_1C$, and the injective
envelop of $S_2$ is $e_2C$. One easily sees that any quotient
comodule of $e_1C$ (or $e_2C$) is left quasi-finite. Hence $C$ is
strictly quasi-finite as a left $C$-comodule (cf. \cite[Theorem
3.1]{GNT}). Similarly, $C$ is also strictly quasi-finite as a right
$C$-comodule. Hence $C$ is artinian. The minimal injective
resolution of $S_1$ is
$$\xymatrix@C=1.0cm{
  0 \ar[r] & S_1 \ar[r] &e_1C\ar[r]^{f}& e_2C \ar[r] & 0, }$$ where $f(x)=y^*\cdot x$ for all $x\in
e_1C$, where $y^*\in C^*$ is the linear map sending $y$ to the unit and other pathes to 0. We want to compute the kernel and
cokernel of the map
$\xymatrix@C=1.5cm{\Hom_{C}(C,e_1C)\ar[r]^{\Hom_{C}(C,f)}&\Hom_{C}(C,e_2C).}$
Let $h_C(-,-)$ be the cohom functor. We have the
commutative diagram of morphisms of left $C^*$-modules (cf. \cite[Appendix]{HTVZ}):
$$\xymatrix@C=1.5cm{h_C(e_1C,C)^*\ar[d]^{\cong}\ar[r]^{h_C(f,C)^*}&h_C(e_2C,C)^*\ar[d]^{\cong}\\
\Hom_{C}(C,e_1C)\ar[r]^{\Hom_{C}(C,f)}&\Hom_{C}(C,e_2C).}$$ By \cite[Appendix]{HTVZ}, we have a commutative diagram of
morphisms of right $C$-comodules:
$$\xymatrix@C=1.5cm{h_C(e_2C,C)\ar[d]^{\cong}\ar[r]^{h_C(y^*\cdot,C)}&h_C(e_1C,C)\ar[d]^{\cong}\\
Ce_2\ar[r]^{\cdot y^*}&Ce_1.}$$
It is clear that the bottom map is
surjective and the kernel of the bottom map is $S_2$. Let $T_1$ and
$T_2$ be the simple right $C$-comodules corresponding to the
vertices. Then, $\Ext_C^0(C,S_1)=0$ and $\Ext^1_C(C,S_1)\cong T_2$
as right $C$-comodules.

Similarly, we see $\Ext_C^0(C,S_2)=0$ and $\Ext^1_C(C,S_2)\cong
T_1$. Hence $C$ is an AS-regular coalgebra, and $C^*$ is an
AS-regular algebra.}
\end{exa}

\section{Calabi-Yau property of AS-regular (co)algebras}

In this section we give the relations
between the CY property and AS-regularity of noetherian complete
semiperfect algebras (or equivalently, artinian coalgebras).

Let $\mathcal{C}$ be a Hom-finite $\k$-linear category. Recall that
an additive functor $\mathscr{S}:\mathcal{C}\to\mathcal{C}$ is call
a {\it right Serre functor} (cf. \cite[Appendix]{Boc}) if there are
natural isomorphisms
$$\eta_{X,Y}:\Hom_{\mathcal{C}}(X,Y)\longrightarrow\Hom_{\mathcal{C}}(Y,\mathscr{S}X),$$
for any $X,Y\in \mathcal{C}$. Note that a Hom-finite $\k$-linear triangulated category is CY-$n$ if and only if the $n$-th shift functor $[n]$ is a right Serre functor.

Let $A$ be a noetherian algebra. Then the triangulated category
$\D_{fd}^b({}_A\M)$ is a Hom-finite $\k$-linear category. If $A$
satisfies the left $\chi$-condition and ${}_{A}A$ has finite injective
dimension, then $\D_{fd}^b({}_A\M)$ has a Serre functor.

\begin{lem} Let $A$ be a noetherian semiperfect algebra with cofinite Jacobson radical. If
$A$ satisfies the left $\chi$-condition and ${}_{A}A$ has finite
injective dimension, then $\mathscr{S}=\RHom_A(-,A)^*$ is a right
Serre functor of $\D_{fd}^b({}_A\M)$.
\end{lem}
\proof We only point out that the $\chi$-condition ensures that
$\mathscr{S}$ is a well defined functor from $\D_{fd}^b({}_A\M)$ to itself. \qed

Let $C$ be an aritinian coalgebra satisfying the left $\chi$-condition.
In this case, we have a natural isomorphism
$\Hom_{C^*}(X,C)\cong X^*$ for every $X\in\D^b_{fd}({}_{C^*}\M)$ in
$\D^b_{fd}(\M_{C^*})$. Indeed, by \cite[Lemma 2.6]{WZ}, we have
$X\cong R\Gamma(X)$ in $\D^b_{fd}({}_{C^*}\M)$. Since $R\Gamma(X)$
is also an object in $\D^b_{fd}(\M^C)$, we have
$\Hom_{C^*}(R\Gamma(X),C)\cong R\Gamma(X)^*$.

\begin{thm}\label{cthem4} Let $C$ be a basic AS-regular coalgebra of global
dimension $n$, and let $\sigma$ be the Nakayama automorphism of $C$.
Then there are natural isomorphisms $$\mathscr{S}(X)\cong
{}_{\sigma^*}X[n],$$ for all $X\in\D^b_{fd}({}_{C^*}\M)$.
\end{thm}
\proof Set $A:=C^*$ and $\alpha:=\sigma^*$. For $X\in
\D^b_{fd}({}_{A}\M)$, we have natural isomorphisms in
$\D^b_{fd}(\M_{A})$:
$$\begin{array}{ccl}
   \mathscr{S}(X)^*&=& \RHom_{A}(X,A)\\
   &\cong&\RHom_{A}(X,{}_1A_{\alpha^{-1}}[n])\ot_A {}_1A_\alpha[-n]\\
   &\overset{(a)}\cong& \Hom_A(R\Gamma(X),C)\ot_A {}_1A_\alpha[-n]\\
   &\cong&R\Gamma(X)^*\ot_A {}_1A_\alpha[-n]\\
   &\cong&X^*\ot_A {}_1A_\alpha[-n].
  \end{array}
$$
The isomorphism $(a)$ above follows from the local duality theorem
\cite[Prop. 3.4]{CWZ} since ${}_1A_{\alpha^{-1}}[n]$ is the balanced
dualizing complex of $A$ (cf. Remark \ref{crem1}). Now by taking the
$\k$-linear dual of the above isomorphisms, we obtain a natural
isomorphism $\mathscr{S}(X)\cong {}_\alpha X[n]$. \qed

Recall that an automorphism of $C$ is {\it inner} if its dual is an
inner automorphism of $C^*$. Clearly, if the Nakayama automorphism
$\sigma$ is inner, then $R\Gamma(C^*)\cong C[-n]$. The converse is
also true. So, immediately, we obtain the following criterion for
the dual algebra of an AS-regular coalgebra to be CY.

\begin{cor}\label{ccor3} Let $C$ be a basic AS-regular coalgebra. If the Nakayama automorphism of $C$ is inner,
then $C^*$ is CY.
\end{cor}

\vspace{5mm}
\subsection*{Acknowledgement} {\small The work is
supported in part by an FWO-grant and NSFC (No. 10801099). The first
named author wishes to thank Department of Algebra and Analysis of
University of Almeria for its hospitality during his visit.}

\bibliography{}

\end{document}